\documentclass[nospthms,envcountsect]{svmult} 
\usepackage{makeidx}     
\makeindex
\usepackage{latexsym}
\usepackage{amsfonts}
\usepackage{amsmath}
\usepackage{amssymb}
\usepackage{amscd}
\usepackage{multicol}    
\usepackage{enumerate}

    \newtheorem{theorem}                    {Theorem}       [section]
    \newtheorem{lemma}      [theorem]       {Lemma}
    \newtheorem{corollary}  [theorem]       {Corollary}
    \newtheorem{proposition}[theorem]       {Proposition}
    \newtheorem{definition} [theorem]       {Definition}
    \newtheorem{conjecture} [theorem]       {Conjecture}

\begin{document}
\catcode`@=11
\atdef@ I#1I#2I{\CD@check{I..I..I}{\llap{$\m@th\vcenter{\hbox
  {$\scriptstyle#1$}}$}
  \rlap{$\m@th\vcenter{\hbox{$\scriptstyle#2$}}$}&&}}
\atdef@ E#1E#2E{\ampersand@
  \ifCD@ \global\bigaw@\minCDarrowwidth \else \global\bigaw@\minaw@ \fi
  \setboxz@h{$\m@th\scriptstyle\;\;{#1}\;$}%
  \ifdim\wdz@>\bigaw@ \global\bigaw@\wdz@ \fi
  \@ifnotempty{#2}{\setbox@ne\hbox{$\m@th\scriptstyle\;\;{#2}\;$}%
    \ifdim\wd@ne>\bigaw@ \global\bigaw@\wd@ne \fi}%
  \ifCD@\enskip\fi
    \mathrel{\mathop{\hbox to\bigaw@{}}%
      \limits^{#1}\@ifnotempty{#2}{_{#2}}}%
  \ifCD@\enskip\fi \ampersand@}
\catcode`@=\active

\renewcommand{\labelenumi}{\roman{enumi})}
\newcommand{\st}{\scriptstyle}
\newcommand{\isom}{\stackrel{\sim}{\longrightarrow}}
\newcommand{\Aut}{\operatorname{Aut}}
\newcommand{\Hom}{\operatorname{Hom}}
\newcommand{\End}{\operatorname{End}}
\newcommand{\Ext}{\operatorname{Ext}}
\newcommand{\Gal}{\operatorname{Gal}}
\newcommand{\Pic}{\operatorname{Pic}}
\newcommand{\Spec}{\operatorname{Spec}}
\newcommand{\trdeg}{\operatorname{trdeg}}
\newcommand{\holim}{\operatornamewithlimits{holim}}
\newcommand{\im}{\operatorname{im}}
\newcommand{\coim}{\operatorname{coim}}
\newcommand{\coker}{\operatorname{coker}}
\newcommand{\gr}{\operatorname{gr}}
\newcommand{\id}{\operatorname{id}}
\newcommand{\Br}{\operatorname{Br}}
\newcommand{\cd}{\operatorname{cd}}
\newcommand{\CH}{\operatorname{CH}}
\renewcommand{\lim}{\operatornamewithlimits{lim}}
\newcommand{\colim}{\operatornamewithlimits{colim}}
\newcommand{\rk}{\operatorname{rank}}
\newcommand{\codim}{\operatorname{codim}}
\newcommand{\NS}{\operatorname{NS}}
\newcommand{\cone}{{\rm cone}}
\newcommand{\rank}{\operatorname{rank}}
\newcommand{\ord}{{\rm ord}}
\newcommand{\f}{{\cal F}}
\newcommand{\p}{{\cal P}}
\newcommand{\N}{{\mathbb N}}
\newcommand{\A}{{\mathbb A}}
\newcommand{\Z}{{{\mathbb Z}}}
\newcommand{\Q}{{{\mathbb Q}}}
\newcommand{\R}{{{\mathbb R}}}
\renewcommand{\H}{{{\mathbb H}}}
\renewcommand{\P}{{{\mathbb P}}}
\newcommand{\F}{{{\mathbb F}}}
\newcommand{\m}{{\mathfrak m}}
\newcommand{\Sm}{{\text{\rm Sm}}}
\newcommand{\Sch}{{\text{\rm Sch}}}
\newcommand{\et}{{\text{\rm et}}}
\newcommand{\eh}{{\text{\rm eh}}}
\newcommand{\Wh}{{\text{\rm Wh}}}
\newcommand{\tra}{{\text{\rm tr}}}
\newcommand{\Zar}{{\text{\rm Zar}}}
\newcommand{\Nis}{{\text{\rm Nis}}}
\newcommand{\tr}{\operatorname{tr}}
\newcommand{\PreShv}{\text{\rm PreShv}}
\newcommand{\Div}{\operatorname{Div}}
\renewcommand{\div}{\operatorname{div}}
\newcommand{\corank}{\operatorname{corank}}
\renewcommand{\O}{{\cal O}}
\renewcommand{\p}{{\mathfrak p}}
\newcommand{\proof}{\noindent{\it Proof. }}
\newcommand{\proofend}{\hfill $\Box$ \\}
\newcommand{\rem}{\noindent {\it Remark. }}
\newcommand{\example}{\noindent {\bf Example. }}
\newcommand{\ar}{{\text{\rm ar}}}

\title*{Arithmetic cohomology over finite fields and special values of
$\zeta$-functions}
\author{Thomas Geisser\thanks{Supported in part by NSF, and the Alfred P. Sloan Foundation}}
\institute{University of Southern California\\
Department of Mathematics, KAP 108\\
3620 Vermont Avenue\\
Los Angeles, CA  90089-2532}

\maketitle

\begin{abstract}
We construct cohomology groups with compact support $H^i_c(X_\ar,\Z(n))$
for separated schemes of finite type over a finite field, which generalize
Lichtenbaum's Weil-etale cohomology groups for
smooth and projective schemes. In particular, if Tate's conjecture
holds, and rational and numerical equivalence agree up to torsion,
then the groups $H^i_c(X_\ar,\Z(n))$ are finitely generated,
form an integral model of $l$-adic cohomology with compact support,
and admit a formula for the special values of the $\zeta$-function
of $X$.
\end{abstract}


\section{Introduction}
In \cite{licht}, Lichtenbaum introduced the Weil-etale topology
in order to produce finitely generated cohomology groups for
varieties over finite fields which are related to special values
of zeta-functions. In \cite{ichweil}, we calculated the precise
relationship between Weil-etale cohomology groups and etale
cohomology groups. In particular, if one assumes Tate's conjecture
on the bijectivity of the cycle map
and Beilinson's conjecture that rational and numerical equivalence
agree up to torsion, then for smooth and projective varieties, the
Weil-etale cohomology groups of the motivic complex have all
properties expected by Lichtenbaum, and allow a new interpretation
of results of Kahn \cite{kahnold}. However, for non-smooth or
non-proper schemes, the Weil-etale cohomology groups are not
finitely generated in general.
In this paper, we use ideas of Voevodsky to construct a
modified version $H^{i}_c(X_\ar ,\Z(n))$ of Weil-etale
cohomology which we call arithmetic cohomology (with compact support).
Arithmetic cohomology groups are expected to be finitely generated, and
related to special values of $\zeta$-functions for every separated scheme
of finite type over a finite field.

To construct arithmetic cohomology groups, we first define an
intermediate Grothendieck topology, called \eh-topology, which is generated by
etale covers and abstract blow-ups. The \eh-topology bears the same
relationship to etale cohomology as Voevodsky's cdh-topology
to the Nisnevich topology. The first advantage of the \eh-topology is that
cohomology groups with compact support
$H^i_c(X_\eh,\f)$ can be defined independently of the choice of a
compactification. An important tool to calculate \eh-cohomology
groups is the following

\begin{theorem} a) For a separated scheme
of finite type $X$ over a field of characteristic not dividing $m$,
$H^i(X_\et,\mu_m^{\otimes n})\cong H^i(X_\eh,\mu_m^{\otimes n})$.

b) Under resolution of singularities for schemes up to dimension $d$,
the etale- and eh-hypercohomology of the motivic complex agree for every
smooth scheme $X$ of dimension at most $d$.
\end{theorem}

The $\eh$-cohomology groups can be used to give a generalization of Tate's
conjecture.

\begin{conjecture}\label{ta}
For every separated scheme of finite type $X$ over $\F_q$, and every $n\in \Z$,
the cycle map from eh-motivic cohomology to the Galois fixed part
of $l$-adic cohomology of $\bar X=X\times_{\F_q}\bar \F_q$,
$$H^i_c(X_\eh,\Z(n))\otimes \Q_l\to
H^i_c(\bar X_\et,\Q_l(n))^{\Gal(\bar \F_q/\F_q)}$$
is an isomorphism, and the Galois-module
$H^i_c(\bar X_\et,\Q_l(n))$ is semi-simple at the eigenvalue $1$.
\end{conjecture}

A homological version has been considered by Jannsen \cite{jannsenln}.
If resolution of singularities and the aforementioned conjectures
of Tate and Beilinson (involving only smooth and projective
schemes) hold, then Conjecture \ref{ta} holds.

Arithmetic cohomology (with compact support) for separated schemes
of finite type over the finite field $\F_q$ is defined by applying
Lichtenbaum's idea
\cite{licht} of replacing the Galois group $\Gal(\bar \F_q/\F_q)$
by the Weil-group $G$ to $\eh$-cohomology. More precisely, a Weil-eh-sheaf 
on $X$ is an eh-sheaf on $\bar X$ together with an action of
$G$, and Weil-eh cohomology (with compact support) $H^{i}_c(X_{\Wh} ,\f)$
of the sheaf $\f$ is defined as the cohomology of the complex
$R\Gamma(G,R\Gamma_c(\bar X_\eh,\f))$.
All results on the relationship between Weil-etale cohomology and etale
cohomology proved in \cite{ichweil} carry over to the present situation.
The groups $H^{i}_c(X_\ar ,\Z(n))$ are defined as the Weil-eh cohomology
groups with compact support of the motivic complex $\Z(n)$ of Suslin-Voevodsky
(we set $\Z(n)=\colim_{p\not| m}\mu_m^{\otimes n}[-1]$ for $n<0$).
Arithmetic cohomology groups are expected to satisfy the following

\begin{conjecture}\label{fg}
For all $n\in \Z$ and schemes $X$ separated and of finite type over $\F_q$,
the groups $H^{i}_c(X_\ar ,\Z(n))$ are finitely
generated, vanish for almost all $i$, and form an integral model for
$l$-adic cohomology with compact support for all $l\not= p$,
$$H^{i}_c(X_\ar ,\Z(n))\otimes \Z_l \cong H^i_c(X_\et,\Z_l(n)). $$
\end{conjecture}

For $l=p$, we get a new theory which agrees with logarithmic de Rham-Witt
cohomology for smooth and proper schemes.
Arithmetic cohomology groups should be related to special values of
zeta-functions in the following way \cite{licht}:

\begin{conjecture}\label{zz}
The weighted alternating sum of the ranks equals the order of the
zeta-function
$$\rho_n:=\sum_i (-1)^i i\cdot \rk H^i_c(X_\ar ,\Z(n))=\ord_{s=n}\zeta(X,s),$$
and for $s\mapsto n$,
$$ \zeta(X,s)\sim\pm
(1-q^{n-s})^{\rho_n}\cdot \chi(H^{*}_c(X_\ar ,\Z(n)),e)
\cdot q^{\chi(n)}.$$
\end{conjecture}

Here $e\in H^{1}((\F_q)_\ar ,\Z(0))\cong \Z$ is a generator,
$\chi(H^{*}_c(X_\ar ,\Z(n)),e)$ is the Euler-characteristic of the complex
$$\cdots \to H^{i-1}_c(X_\ar ,\Z(n))\stackrel{\cup e}{\longrightarrow}
H^{i}_c(X_\ar ,\Z(n))\stackrel{\cup e}{\longrightarrow}
H^{i+1}_c(X_\ar ,\Z(n))\to\cdots,$$
and $\chi(n)=
\sum_{0\leq i\leq n\atop j}(-1)^{i+j}(n-i)\cdot \dim H^j_c(X_\eh,\Omega^i)$
is a generalization of the correcting factor of Milne \cite{milnevalues}.

We show several implications between the above conjectures, and
the following

\begin{theorem} Conjectures \ref{ta}, \ref{fg} and \ref{zz} hold
in the following instances:

a) For every $n$, if $X$ is a curve.

b) For $n\leq 0$, if $X$ is a scheme of dimension at most $d$, and
resolution of singularities for schemes of dimension at
most $d$ exists.

c) For every $X$ and every $n$, if the conjectures of Tate and Beilinson
hold, and resolution of singularities exists.
\end{theorem}

Finally, we give an example showing that
essentially every statement given above is incorrect if one
uses the etale topology.

\medskip

{\it Acknowledgments}: The results of this paper were obtained
during a stay at the University of Tokyo. I would like to express
my gratitude to the University of Tokyo, and especially to
T.Saito for their hospitality and financial support.
I thank C.Weibel for help with the example in
Proposition \ref{counterex}, to B.Kahn for comments on an
earlier version of the paper, and to the referee, whose careful
reading improved the exposition.

\section{The etale h-topology}
We introduce a Grothendieck topology which is finer than
the etale topology, and has several advantages over
the etale topology. For instance, there are well-defined cohomology
groups with compact support, and there is a
long exact sequence of cohomology groups for blow-ups.
On the other hand, for locally constant
torsion coefficients, one gets the same cohomology groups as
for the etale topology.

We use the term Grothendieck topology on a subcategory $\cal C$ of the
category of schemes in the sense of \cite[II 1]{sga4}.
In particular, any morphism $Y\to X$ that can be dominated
by a covering $U\to X$ is itself a covering. If we are given
for every object $X$ of $\cal C$ a class of morphisms with target
$X$, then the intersection of all topologies containing these
morphisms is again a topology, called the Grothendieck topology
generated by these morphisms \cite[II 1.1.6]{sga4}.

\begin{definition}
The etale h-topology (or short eh-topology) on a suitable
subcategory of the category of schemes is the Grothendieck topology
generated by the following coverings:

1) etale coverings

2) abstract blow-ups: Assume we have a cartesian square
\begin{equation}\begin{CD}\label{blowup}
Z'@>i'>> X'\\
@Vf'VV @VfVV \\
Z@>i>> X,\end{CD}\end{equation}
where $f$ is proper, $i$ is a closed embedding,
and $f$ induces an isomorphism $X'-Z'\stackrel{\sim}{\longrightarrow} X-Z$.
Then $(X'\stackrel{f}{\to} X,Z \stackrel{i}{\to} X)$ is a covering.
\end{definition}

The definition is motivated by Voevodsky's cdh-topology,
which is generated by Nisnevich covers and abstract blow-ups.
For singular schemes, the cdh-topology has better properties than the
Nisnevich topology, and similarly the
\eh-topology has better properties than the etale topology.

\smallskip

\noindent{\it Example.} Every scheme $X$ is covered by its irreducible
components, $X^{red}\to X$ is a covering, and for every
blow up $X'$ of $X$ with center $Z$, $(X'\to X,Z \to X)$
is a covering.
\smallskip

\rem It is tempting to use the h-topology of Voevodsky instead of
the \eh-topology, in order to use alterations of de Jong instead of
resolution of singularities.
However, if one does so, then one looses the mod $p$ information.
Similarly, Voevodsky's method to define motivic cohomology groups
as Ext-groups in the derived category of mixed motives does not
give well-behaved cohomology groups for the etale topology
(the resulting cohomology groups will be $p$-divisible).
Even with rational coefficients we do not know how to prove the
analog of Proposition \ref{blowuples} for an alteration, and hence
we cannot prove Theorem \ref{sameasbefore} for finer topologies
than the eh-topology.

\smallskip

We recall the following facts from \cite{sv00}:

\begin{lemma}\label{refines}
a) Every proper morphism $p:X'\to X$, such that for every point $x\in X$
there is a point in $p^{-1}(x)$ with the same residue field, is an
eh-covering. 

b) Every abstract blow-up $(X'\to X, Z\to X)$ has a refinement
$(\tilde X\to X, Z\to X)$ such that every irreducible component of $\tilde X$
has dimension not larger than the dimension of $X$.
\end{lemma}

\proof
a) \cite[Lemma 5.7]{sv00}

b) We can replace $X'$ by the disjoint union of its
irreducible components. We claim that removing all irreducible components
of dimension larger than the dimension of $X$ from $X'$ gives a
refinement $(\tilde X\to X'\to X, Z\to X)$. Indeed, let $T$ be one of the
irreducible components, and let $\eta$ be the generic point of $T$.
If $\eta$ maps to $X-Z$ under $f$, then $T$ is birational to a component to
$X$, hence $\dim T\leq \dim X$.
If $\eta$ maps to $Z$, then the map $T\to X$ factors through $Z$,
and we can remove $T$, because by a) the resulting scheme is still a
covering.
\proofend

A covering as in a) is called proper eh-covering.
If $X$ is integral, then a proper \eh-covering which is an isomorphism
over a neighborhood of the generic point of $X$ is called a
proper birational \eh-covering. Every proper \eh-cover of an integral
scheme $X$ admits a proper birational refinement by the argument
in \cite[Lemma 5.7]{sv00}.

\begin{proposition}\label{normalformgeneral}
Every eh-cover of $X$ has a refinement of the form
$$ \{ U_i\to X'\to X \}_{i\in I},$$
where $\{U_i\to X'\}_{i\in I}$ is an etale cover, and $X'\to X$
is a proper eh-cover.
\end{proposition}

\proof \cite[Prop. 5.9]{sv00}
\proofend

We now fix a perfect field $k$, let $\Sch/k$ be the category
of separated schemes of finite type over $k$, and $\Sm/k$ the
full subcategory of smooth schemes.
For $d\in {\mathbb N}\cup {\infty}$,
we denote by $\Sch^d/k$ and $\Sm^d/k$ the full subcategory consisting
of schemes of Krull dimension at most $d$ of $\Sch/k$
and $\Sm/k$, respectively. By Lemma \ref{refines} b), we can consider the
\eh-topology on $\Sch^d/k$. We write $(\Sch^d/k)_{\eh}$ for the
category $\Sch^d/k$ equipped with the \eh-topology, and
$(\Sch^d/k)_\eh^\sim$ for the topos of \eh-sheaves on $\Sch^d/k$.
Similarly, $(\Sm^d/k)_\et$ is the category of smooth schemes equipped
with the etale topology, and $(\Sm^d/k)_\et^\sim$ is the topos of
etale sheaves on $\Sm^d/k$.

\begin{definition}
For $d\in {\mathbb N}\cup \infty$ we denote by $R(d)$
the strong form of resolution of singularites for varieties
up to dimension $d$, i.e. the following two conditions:

\begin{itemize}
\item
For every integral separated scheme $X\in \Sch^d/k$,
there is a proper, birational map $f:Y\to X$ with $Y\in \Sm/k$.

\item
For every smooth scheme $X\in \Sm^d/k$ and
every proper birational map $f:Y\to X$, there is a sequence
of blow-ups along smooth centers $X_n\to X_{n-1}\to \cdots \to X_1\to X$
such that the composition $X_n\to X$ factors through $f$.
\end{itemize}
\end{definition}

If $\text{char}\;k=0$, $R(\infty)$ holds
by Hironaka's theorem. By Abhyankar, $R(2)$ is known
in general, and $R(3)$ is known for algebraically closed fields of
characteristic $p>5$.
Condition $R(d)$ implies that every scheme in $\Sch^d/k$
is locally smooth for the \eh-topology.

Let $\rho_d:(\Sch^d/k)_\eh \to (\Sm^d/k)_\et$ be the canonical maps
of sites, and for $a\geq b$ let 
$\iota:(\Sch^a/k)_\eh^\sim \to (\Sch^b/k)_\eh^\sim$
and $\sigma:(\Sm^a/k)_\eh^\sim \to (\Sm^b/k)_\eh^\sim$
be the canonical morphism of topoi induced by the restriction map.

\begin{lemma}\label{topoilemma} Assume that $R(d)$ holds.

a) The functor $\rho_d$ induces a morphism of topoi 
$\rho_d:(\Sch^d/k)_\eh^\sim \to (\Sm^d/k)_\et^\sim$.

b) For every $a\geq b$, there is a commutative diagram
$$\begin{CD}
(\Sch^a/k)_\eh^\sim @<\rho_a^*<< (\Sm^a/k)_\et^\sim\\
@V\iota_*VV @V\sigma_* VV\\
(\Sch^b/k)_\eh^\sim @<\rho_b^*<< (\Sm^b/k)_\et^\sim.
\end{CD}$$
\end{lemma}

\proof
a) The only point which needs explanation is the left exactness of $\rho_d^*$.
On a smooth scheme $S$, the presheaf pull-back $\rho^p_d\f(S)$ agrees with 
$\f(S)$.
Since the system of coverings occuring in the definition of the 
sheafification functor can be assumed to be filtered by 
\cite[III 2.2 a)]{milnebook}, and since coverings by smooth schemes 
are cofinal in the system of all eh-covers by $R(d)$, the statements follows
\cite[App. A, Prop. 7]{milnebook}

b) By the same argument as in a), it suffices to know the value of the 
presheaf pull-back on smooth schemes of dimension at 
most $a$ in order to calculate $\rho_b^*$ and $\rho_a^*$. But on such a 
scheme,  $\rho_b^p\sigma_*\f=\f=\iota_*\rho_a^p\f$.
\proofend

We do not know that the presheaf pull-back $\rho^p_d$ is left exact.
This is because equalizers do not exist in the category of smooth schemes, 
so that
the colimit system defining the presheaf pull-back is not filtered, see
\cite[Rem. 3.7]{friedvoe}.

By resolution of singularities, every proper birational map
to a smooth scheme can be refined by blow-ups along smooth
centers. Since a blow-up $X'\to X$ along a smooth center
satisfies the hypothesis of Lemma \ref{refines} a), Proposition
\ref{normalformgeneral} gives

\begin{corollary}\label{normalform}
Let $X\in \Sm^d/k$ and assume that condition $R(d)$ holds.
Then every eh-cover of $X$ has a refinement of the form
$$ \{ U_i\to X'\to X \},$$
where $\{U_i\to X'\}$ is an etale cover, and $X'\to X$
is a composition of blow-ups along smooth centers.
\end{corollary}

The following lemma will be applied several times.

\begin{lemma}\label{devissage}
(Devissage Lemma) Let $P(X)$ be a property for schemes in $\Sch^d/k$. Assume
the following:
\begin{enumerate}
\item Condition $R(d)$.
\item If $Z\subseteq X$ is a closed subscheme of $X$ with
open complement $U$, and if $P(-)$ holds for two of the
three schemes $X,U$ and $Z$, then it also holds for the third.
\item $P(X)$ holds for all smooth and projective $X\in \Sm^d/k$.
\end{enumerate}
Then $P(X)$ holds for all $X\in \Sch^d/k$.
\end{lemma}

\proof We proceed by induction on the dimension of $X$.
Given $X$, we can assume by noetherian
induction that $P(Z)$ holds for all closed subschemes $Z$ of $X$,
and using property ii) we can reduce to the case that $X$ is integral.
By Chow's lemma, there is a
projective scheme $X_1$ and an open subscheme $U_1$ of $X_1$
isomorphic to an open subscheme of $X$. Let $X_2$ be a
desingularization of $X_1$, then there is an open subscheme
$U_2$ of $X_2$ isomorphic to an open subscheme of $X$.
By condition iii) we have $P(X_2)$, and by condition ii)
this implies $P(U_2)$ and then $P(X)$.
\proofend

\section{Cohomology for the eh-topology}

The usual argument with generators \cite[III Lemma 1.3]{milnebook}
shows that the categories $(\Sch^d/k)_\eh^\sim$ have enough injectives.
Given an eh-sheaf $\f\in (\Sch^d/k)_\eh^\sim$ and $X\in \Sch^d/k$,
the cohomology groups $H^i(X_\eh,\f)$ are defined
as the derived functors of the global section functor
$\f\to \Gamma(X_\eh,\f)$.
For $X\in \Sch/k$, we let $\Z(X)$ be the free presheaf
$U\mapsto \Z[\Hom_\Sch(U,X)]$ represented by the scheme $X$,
and let $\Z_\eh(X)$ be its associated \eh-sheaf.
Sheafification is necessary, because the \eh-topology
is not subcanonical.

\begin{lemma}\label{poiu} For every $a\geq b$, the canonical morphism of topoi
$\iota:(\Sch^a/k)_\eh^\sim \to (\Sch^b/k)_\eh^\sim$ 
induces an isomorphism of cohomology groups. Moreover,
$H^i(X_\eh,\f)=\Ext^i_{(\Sch^d/k)_\eh^\sim}(\Z_\eh(X),\f)$
for every eh-sheaf $\f\in (\Sch^d/k)_\eh^\sim$
of abelian groups and $X\in \Sch^d/k$.
\end{lemma}

\proof
The first statement is proved as in \cite[III Prop. 3.1]{milnebook}:
Clearly $\iota_*$ is exact, and it suffices to show that
$\f\cong \iota_*\iota^*\f$. This is clear on the presheaf level,
by definition of the presheaf pull-back, and $\iota_*\iota^*\f$ 
is a sheaf on $(\Sch^d/k)_\eh$ by Lemma \ref{refines} b).
The second statement follows because
$\Hom_{(\Sch^d/k)_\eh^\sim}(\Z_\eh(X),\f)\cong \f(X)$
for every \eh-sheaf $\f\in (\Sch^d/k)_\eh^\sim$
of abelian groups and $X\in \Sch^d/k$.
\proofend

\begin{proposition}\label{lexprop}
Every abstract blow-up square \eqref{blowup}, gives rise to
a long exact sequence of cohomology groups:
\begin{equation}\label{lex}
\cdots \to H^i(X_\eh,\f)\stackrel{i^*,f^*}{\longrightarrow}
H^i(Z_\eh,\f)\oplus H^i(X'_\eh,\f)\stackrel{{f'}^*-{i'}^*}{\longrightarrow}
H^i(Z'_\eh,\f)\to \cdots.
\end{equation}
\end{proposition}

\proof
By Lemma \ref{poiu}, it suffices to show that there is a short exact sequence
of \eh-sheaves
\begin{equation}\label{referee}
0\to \Z_\eh(Z')\stackrel{f'_*,i'_*}{\longrightarrow}
\Z_\eh(Z)\oplus \Z_\eh(X') \stackrel{i_*-f_*}{\longrightarrow}\Z_\eh(X)\to 0.
\end{equation}
Exactness on the clear on the presheaf level because
$i'$ is injective. To check exactness in the middle, let
$x\in \Z_\eh(Z)(U)\oplus \Z_\eh(X')(U)$ be a section of the middle term
over a scheme $U$. By going to an \eh-cover, we can assume that $U$ is
integral and that $x=(\sum_ln_l\alpha_l,\sum_jm_j\beta_j)$
is represented by a linear combination of pairwise different morphisms
$\alpha_l:U\to Z$ and pairwise different morphisms $\beta_j:U\to X'$
such that $\sum_ln_li\circ \alpha_l=\sum_jm_jf\circ \beta_j\in \Z(X)(U)$.
If $\beta_j$ maps the generic point of $U$ to $Z'$, then $\beta_j$
factors through $Z'$, hence changing $x$ by an element in the image
of $\Z(Z')(U)$, we can assume that every summand $\beta_j$ of $x$ sends
the generic point of $U$ to $X'-Z'$. We claim that this implies
that $x=0$, because no two $\alpha_l$ or $\beta_j$ can become
equal in $\Z(X)(U)$. This is clear for the $\alpha_l$
since $i:Z\to X$ is injective and the $\beta_j$ don't have image
in $Z$. On the other hand, if $\beta_1$ and $\beta_2$ are two maps
which map the generic point of $U$ to $X'-Z'$, and which become equal
when composed with $f$, then because
$f:X'-Z'\to X-Z$ is an isomorphism, $\beta_1$ and $\beta_2$ agree on
the generic point of $U$, hence are equal.

To show exactness on the right, let $f\in \Z(X)(U)$ be a morphism
and consider the eh-covering $U\times_XZ,U\times_X{X'}$ of $U$.
Restricting $f$ to this covering, we get two maps 
$f_Z:U\times_XZ\to U\to X$ and $f_{X'}:U\times_X{X'}\to U \to X$, 
which clearly lie in the image of $\Z(Z)(U\times_XZ)$ and
$\Z(X')(U\times_XX')$, respectively.
\proofend

\rem
The presheaf analog of \eqref{referee} is not exact in the middle,
as stated in \cite[Lemma 12.1]{sv00}. Take for example $U=\Spec k[x,y]/(xy)$,
$X$ the affine plane and $X'$ be the blow-up of $X$ at the origin.
Take $g$ to be the map $U\to X'$ which embeds the line $x=0$ into the
exceptional divisor $Z'$, and maps the line $y=0$ to any line of $X'$
intersecting the exceptional divisor in the image of the origin
$(0,0)\in U$.
If $\tau$ is the reflection of $U$ sending $y$ to $-y$, then
$g-g\circ \tau$ becomes zero when composed with
the projection $X'\to X$, but does not factor through $Z'$.

We now come to the definition of cohomology with compact support.

\begin{definition}\label{definecompact}
Let $U\in \Sch^d/k$ and let $j:U\to X$ be a compactification of $U$,
i.e. a proper scheme over $k$ containing $U$ as a dense open
subscheme. Let $i:Z\to X$ be the closed complement of $U$ with the reduced 
subscheme structure. For an eh-sheaf $\f\in (\Sch^d/k)_\eh^\sim$,
let $\f \to I^\cdot$ be a resolution by injective eh-sheaves. We
define the eh-cohomology with compact support of $U$ to be
$$R\Gamma_c(U_\eh,\f)=\cone(I^\cdot(X)\to I^\cdot(Z))[-1].$$
\end{definition}

\begin{lemma}
The above definition is independent of the choice of $X$.
\end{lemma}

\proof
Given two compactifications $X$ and $X'$, we can by the usual argument assume 
that there is a map $f:X'\to X$ which is the identity on $U$.
Let $Z'=X'-U\stackrel{i'}{\to} X'$ and $Z=X-U\stackrel{i}{\to}X$ be the
closed complements of $U$ in $X'$ and $X$, respectively, with the
reduced subscheme structure. Since $Z'\cong Z\times_XX'$ as topological spaces,
and since the \eh-cohomology of $Z'$ and $(Z')^{\text{red}}$ agree,
we can assume that $Z'\cong Z\times_XX'$, so that $Z',X',Z$ and $X$ form a
blow-up square \eqref{blowup}. Consider the diagram
$$\begin{CD}
I^\cdot(X)@>i^*>> I^\cdot(Z)@>>> \cone(i^*)\\
@Vf^*VV @V{f'}^*VV @VVV\\
I^\cdot(X') @>{i'}^*>> I^\cdot(Z')@>>> \cone({i'}^*).
\end{CD}$$
Applying $\Hom(-,I^\cdot)$ to the exact sequence \eqref{referee},
we see that the right vertical map is an isomorphism.
\proofend

As usual, \eh-cohomology groups with compact support are contravariant
for proper maps and covariant for open embeddings.
For a closed embedding
$Z\subseteq X$ with open complement $U$ there is a long exact sequence
\begin{equation}\label{cs}
\cdots \to H^j_c(U_\eh,\f)\to H^j_c(X_\eh,\f)\to H^j_c(Z_\eh,\f)
\to \cdots.
\end{equation}

\rem
Another approach to define cohomology with compact support is to let the
$\Z^c_\eh(X)$ be the eh-sheaf associated to the presheaf which sends an
irreducible scheme $V$ to the free abelian group on closed subschemes
$Z\subseteq V\times X$ such that the projection $Z\to V$ is an open embedding.
If $X$ is proper, then $\Z^c_\eh(X)=\Z_\eh(X)$, because then $Z\cong V$
can be identified with the graph of a morphism $V\to X$.
Hence if one defines cohomology with compact support of the sheaf $\f$
as $H^i_c(X_\eh,\f)=\Ext^i_\eh(\Z^c_\eh(X),\f)$, then this agrees with
Definition \ref{definecompact} by Lemma \ref{poiu}.
For an open subscheme $U\subseteq X$ with closed complement $Z$,
one can prove as in \cite[Prop. 3.8]{friedvoe} that there is a short exact
sequence of \eh-sheaves (this fails for the etale topology)
$$0\to \Z^c_\eh(Z) \to \Z^c_\eh(X)\to \Z^c_\eh(U)\to 0,$$
hence the two definitions agree in general in view of \eqref{cs}.

\begin{lemma}\label{ehsheaf}
Let $\f\in (\Sch/k)_\et^\sim$ be a constructible sheaf of
abelian groups. In the blow-up square \eqref{blowup}, let $g$ be the
composition $Z'\to X$. Then there is an exact triangle in the derived
category of etale sheaves on $X$,
$$ \f\to Rf_*\f\oplus i_*\f\to Rg_*\f \to \f[1].$$
In particular, $\f$ is a sheaf for the eh-topology, and there
is a long exact sequence \eqref{lex} for the etale cohomology
of $\f$.
\end{lemma}

\proof
Since $\f$ is constructible, $\f$ is torsion and $i^*\f$, $f^*\f$,
and $g^*\f$ are the restrictions of $\f$ to $Z$, $X'$, and $Z'$, respectively.
It suffices to show that there is a short exact sequence of etale sheaves
$$ 0\to \f\to f_*\f\oplus i_*\f\to g_*\f \to 0,$$
and isomorphisms $R^sf_*\f\cong R^s g_*\f$ for $s>0$.
If $x$ is a geometric point over $X-Z$, then $(R^sg_*\f)_x=(i_*\f)_x=0$
for all $s\geq 0$. Since $f$ is an isomorphism in a neighborhood of $x$,
$(R^sf_*\f)_x=0$ for $s>0$, and the stalk at $x$ of the sequence becomes
the isomorphism $(\f)_x\isom (f_*\f)_x $. For $x$ a geometric point over
$Z$, there are isomorphisms $(\f)_x\isom(i_*\f)_x $ and
$(R^sf_*\f)_x\isom (R^sg_*\f)_x$ by the proper base-change theorem.

Finally, if $\f$ satisfies the sheaf property
for a class of morphisms, then it also satisfied the sheaf property for
the Grothendieck topology generated by it.
\proofend

If the abstract blow-up $f:X'\to X$ is finite, then in the Lemma
it suffices to assume that $\f$ is locally constructible.
Indeed, in this case $R^sf_*\f= R^s g_*\f=0$ for $s>0$, and
the proper base-change theorem is not needed.

Consider the canonical morphism of topoi
$\tau: (\Sch/k)_\eh^\sim\to (\Sch/k)_\et^\sim$.

\begin{theorem}\label{susvoeequal}
Let $\f\in (\Sch/k)_\et^\sim$ be a constructible sheaf. Then
$R^s\tau_*\f=0$ for $s>0$. In particular, for every $X\in \Sch/k$,
$$H^s(X_\et,\mu_m^{\otimes n})\cong H^s(X_\eh,\mu_m^{\otimes n}).$$
\end{theorem}

\proof
By the Lemma, $\f$ is a sheaf for the eh-topology, and we identify 
$\tau^*\f$ with $\f$. Let $C^\cdot$ be the cone of the canonical map of
complexes of etale sheaves $\f\to R\tau_*\f$. It suffices to show that 
$H^i(X_\et,C^\cdot)=0$ for every scheme $X$ in $\Sch/k$. Assume we have
a non-zero element $0\not= u\in H^i(X_\et,C^\cdot)$,
and that $X$ is a scheme of smallest dimension admitting
such an element. Given an abstract blow-up diagram \eqref{blowup},
then according to Propositions \ref{lexprop} and Lemma \ref{ehsheaf},
there is a map of long exact sequences
$$\begin{CD}
H^i(X_\et,\f)@>>> H^i(Z_\et,\f)\oplus H^i(X'_\et,\f)
@>>>  H^i(Z'_\et,\f) \\
@V\tau_XVV @V\tau_ZV\tau_{X'}V @V\tau_{Z'}VV \\
H^i(X_\eh,\f)@>>>  H^i(Z_\eh,\f)\oplus H^i(X'_\eh,\f)
@>>>  H^i(Z'_\eh,\f)
\end{CD}$$
If $X'$ is an irreducible component of $X$ and $Z$ the union
of the remaining components, and if $\tau_X$ is
not an isomorphism, then either $\tau_{X'}$ or $\tau_Z$ is not an
isomorphism, because $\tau_{Z'}$ is an isomorphism by minimality
of the dimension of $X$. Hence we can by induction on the number of
irreducible components of $X$ assume that $X$ is integral.

Since $\tau^*C^\cdot=0$, there is an eh-covering
of $X$ such that $C^\cdot$ is quasi-isomorphic to zero when restricted to this
covering. We can by Proposition \ref{normalformgeneral} assume that the
covering is a composition of an etale cover $\{U_i\to X' \}$, and a proper
eh-cover $X'\to X$.
Replace $f:X'\to X$ by a proper birational refinement, and let $Z$
be the closed subscheme of $X$ where $f$ is not an isomorphism.
Then we get a diagram as above, and by minimality of $X$, $\tau_Z$
and $\tau_{Z'}$ are isomorphisms, hence $\tau_{X'}$ cannot be
an isomorphism for all $i$, and thus $C^\cdot|_{X'}$
is not quasi-isomorphic to zero. But then it is also not
quasi-isomorphic to zero on the etale cover $\{U_i\to X'\}$,
a contradiction.
\proofend

\section{Motivic, Hodge and de Rham cohomology}
Consider the restriction of the motivic complex 
$\Z(n)\in K^-((\Sm^d/k)_\et)$, 
$n\geq 0$, of Suslin-Voevodsky \cite[Def. 3.1]{sv00}, a bounded above complex
of etale sheaves of abelian groups on $\Sm^d/k$. Abusing notation,
we simply write $\Z(n)$ for the extension 
$\rho_d^*\Z(n)\in K^-((\Sch^d/k)_\eh)$. Under $R(a)$, 
the cohomology of $\rho_d^*\Z(n)$ does not depend on $d$ as long as $d\leq a$
by Lemmas \ref{topoilemma} b) and \ref{poiu}.
For negative $n$, we set $\Z(n)=\colim_{p\not|m}\mu_m^{\otimes n}[-1]$.
For an abelian group $A$, we write $A(n)$ for the complex $A\otimes \Z(n)$. 
If $A$ is torsion free, then $H^i(X_\eh,A(n))\cong 
H^i(X_\eh,\Z(n))\otimes A$.

\begin{lemma}\label{sheafisom}
Under $R(d)$, there are quasi-isomorphism
$\Z(0)\cong \Z$, $\Z(1)\cong {\mathbb G}_m[-1]$ and
$\Z/m(n)\cong \mu_m^{\otimes n}$ for all $n\in \Z$ and
$\text{char}\; k\not| m$ in $D^-((\Sch^d/k)_\eh)$.
In particular, $H^i_c(X_\eh,\Z/m(n))\cong H^i_c(X_\et,\mu_m^{\otimes n})$.
\end{lemma}

\proof
This follows by exactness of $\rho_d^*$ from the corresponding
statement for smooth schemes \cite[Lemma 3.2]{sv00}, because $\rho_d^*\Z=\Z$,
$\rho_d^*{\mathbb G}_m={\mathbb G}_m$ and
$\rho_d^*\mu_m^{\otimes n}=\mu_m^{\otimes n}$.
The final statement follows with Theorem \ref{susvoeequal}.
\proofend

\begin{proposition}\label{blowuples}
Let $n\in \Z$, $X\in \Sm/k$, $X'$ the blow-up of $X$ along the smooth
center $Z$, and $Z'=Z\times_X X'$. Then there is a long exact
sequence
$$ \cdots \to H^i(X_\et,\Z(n))\to  H^i(Z_\et,\Z(n))\oplus H^i(X'_\et,\Z(n))
\to  H^i(Z'_\et,\Z(n))\to \cdots .$$
\end{proposition}

\proof
It suffices to prove exactness
rationally, with mod $p^r$ and with mod $m$-coefficients for
$p\not|m$. For mod $m$-coefficients, we have $\Z/m(n)\cong \mu_m^{\otimes n}$,
hence the claim follows from Lemma \ref{ehsheaf}. The result with mod
$p^r$-coefficients is \cite[IV \S1]{grosthesis} because of
$\Z/p^r(n)\cong W_r\Omega_{X,log}^n$ \cite{marcI}.
Finally, rationally motivic cohomology and etale motivic
cohomology agree \cite[Prop. 3.6]{ichdede}, hence the result with rational
coefficients is \cite[Lemma IV 3.12]{marcbook}.
\proofend

\begin{theorem}\label{sameasbefore}
If $n\in \Z$, $X\in \Sm^d/k$ and condition $R(d)$ holds, then
$$H^i(X_\et,\Z(n))\cong H^i(X_\eh,\Z(n)).$$
\end{theorem}

\proof
Let $C^\cdot$ be the cone of the canonical map of
complexes of etale sheaves $\Z(n)\to R(\rho_d)_*\Z(n)$.
It suffices to show that $H^i(X_\et,C^\cdot)=0$
for every scheme $X\in \Sm^d/k$. Assume we have
a non-zero element $0\not= u\in H^i(X_\et,C^\cdot)$,
and that $X$ is a scheme of smallest dimension admitting
such an element. Since $\rho_d^*C^\cdot=0$, there is an eh-covering
of $X$ such that $C^\cdot$ is quasi-isomorphic to zero when restricted to this
covering. By Corollary \ref{normalform} we can assume that the covering has
is a composition of an etale cover $\{U_i\to X' \}$, and a composition of
blow-ups along smooth centers $X'\to X$. Given a blow-up $X'$ of $X$
along the smooth center $Z$, we can find by Propositions \ref{lexprop}
and \ref{blowuples}, a map of long exact sequences
$$\begin{CD}
H^i(X_\et,\Z(n))@>>> H^i(Z_\et,\Z(n))\oplus H^i(X'_\et,\Z(n))
@>>>  H^i(Z'_\et,\Z(n)) \\
@V\tau_XVV @V\tau_ZV\tau_{X'}V @V\tau_{Z'}VV \\
H^i(X_\eh,\Z(n))@>>>  H^i(Z_\eh,\Z(n))\oplus H^i(X'_\eh,\Z(n))
@>>>  H^i(Z'_\eh,\Z(n)).
\end{CD}$$
By minimality of $X$, $\tau_Z$ and $\tau_{Z'}$ are isomorphisms, and we
conclude that $u|_{X'}$ is non-zero. In particular, $C^\cdot|_{X'}$
is not quasi-isomorphic to zero. But then it is also not
quasi-isomorphic to zero on the etale cover $\{U_i\to X'\}$,
a contradiction.
\proofend

\begin{corollary}\label{finite}
Assume $R(d)$, let $n\in \Z$ and $X\in \Sch^d/k$.
If $\cd_l(k)<\infty$ for all primes $l$ dividing $m$, then the groups
$H^i_c(X_\eh,\Z/m(n))$ are finitely generated.
\end{corollary}

\proof
By Lemma \ref{devissage} we can assume that $X$ is smooth
and projective. Write $m=m'\cdot p^r$ with $p\not|m'$. Finite generation of 
$H^i_c(X_\et,\mu_{m'}^{\otimes n})$ is well-known. 
On the other hand, by Theorem \ref{sameasbefore} and \cite{marcI},
$H^i_c(X_\eh,\Z/p^r(n))\cong H^i(X_\et,\Z/p^r(n))
\cong H^{i-n}(X_\et,\nu_r^n)$. The latter group is finitely generated
by \cite[Prop. 4.18]{grossuwa}.
\proofend

\begin{proposition}
Under resolution of singularities, the rational eh-motivic
cohomology groups $H^i(X_\eh,\Q(n))$ agree
with Voevodsky's rational motivic cohomology groups
$\Hom_{DM^-}(M(X),\Q(n)[i])$.
\end{proposition}

\proof
By \cite[Thm. 1.5]{sv00}, Voevodsky's motivic cohomology groups agree with the
Nisnevich cohomology groups of the motivic complex for smooth $X$,
$\Hom_{DM^-}(M(X),\Z(n)[i])= H^i(X_{\text {Nis}},\Z(n))$ .
Using the fact that rationally Nisnevich and
etale cohomology agree, this implies the claim for smooth $X$:
$$\begin{CD}
H^i(X_{\text{Nis}},\Q(n))@>\sim >> \Hom_{DM^-}(M(X),\Q(n)[i])\\
@V\cong VV \\
H^i(X_\et,\Q(n))@>\sim >> H^i(X_\eh,\Q(n)).
\end{CD}$$
For arbitrary $X$, we can proceed by induction on the dimension of
$X$ from the smooth case, using the blow-up sequences \eqref{lex}
for both theories.
\proofend


\begin{proposition}\label{projb}
Assume $R(d+r)$, let $n\in \Z$ and $X\in Sch^d/k$.

a) (Projective bundle formula)
There are canonical isomorphisms
$$
\bigoplus_{j=0}^r H^{i-2j}_c(X_\eh,\Z(n-j))\cong
H^i_c(({\mathbb P}^r_X)_\eh,\Z(n)).
$$

b) (Affine bundle formula) There are canonical isomorphisms
\begin{equation}\label{homotopy}
H^i_c(({\mathbb A}^r_X)_\eh,\Z(n))\cong H^{i-2r}_c(X_\eh,\Z(n-r)).
\end{equation}
\end{proposition}

\proof
a) Comparing the long exact sequences for cohomology
with compact support \eqref{cs}, we can by Lemma \ref{devissage} assume
that $X$ is smooth and projective. In this case, we can consider
etale cohomology instead of eh-cohomology by Theorem \ref{sameasbefore}.
The isomorphism is given by $\sum p_r^*\cup \xi^j$, where
$p_r:{\mathbb P}^r_X\to X$ is the projection, and
$\xi\in H^2(({\mathbb P}^r_X)_\et,\Z(1))\cong \Pic {\mathbb P}^r_X$ is the 
class of ${\cal O}(1)$.

With mod $m$-coefficients, $p\not| m$, this
map is an isomorphism by \cite[Prop. 10.1]{milnebook} because
$\Z/m(n)_\et\cong \mu_m^{\otimes n}$ for every $n\in \Z$. For
$p^r$-coefficients, the map is an isomorphism by
\cite[I Thm. 2.2.11]{grosthesis} because
$\Z/p^r(n)_\et\cong W_r\Omega^n_{X,log}$.
Rationally, the map is an isomorphism by \cite[Thms. 1.5, 4.5]{sv00},
because motivic cohomology and etale motivic cohomology agree.

b) Again we first reduce to the case that $X$ is smooth and projective.
In this case, by \eqref{cs}, the section
${\mathbb P}^{r-1}_X\stackrel{0}{\longrightarrow} {\mathbb P}^r_X$
gives a long exact sequence
\begin{equation}\label{projlong}
\cdots \to H^i_c(({\mathbb A}^r_X)_\eh,\Z(n)) \to
H^i_c(({\mathbb P}^r_X)_\eh,\Z(n))\stackrel{0^*}{\longrightarrow}
H^i_c(({\mathbb P}^{r-1}_X)_\eh,\Z(n)) \to \cdots.
\end{equation}
It is easy to see that $0^*\xi=\xi$, so that
$0^*(p_r^*\cup\xi^j)=p_{r-1}^*\cup \xi^j$ for $0\leq j\leq r-1$.
Consequently, $0^*$ is surjective, and
$\ker 0^*\cong H^{i-2r}_c(X_\eh,\Z(n-r))$.
\proofend

\rem
Since $H^i(({\mathbb P}_k^1)_\et,\Z/m(0))\cong
H^i(k_\et,\Z/m)\oplus H^i(k_\et,\mu_m^{\otimes -1})$,
the projective bundle formula cannot hold in general without
modifying the definition of motivic cohomology for negative $n$.

\rem
In view of Proposition \ref{projb} b) and in analogy with the
formula for the $\zeta$-function, it would be natural
to define negative eh-cohomology as
$$H^i_c(X_\eh,\Z(n))=H^{i-2n}_c(({\mathbb A}^{-n}_X)_\eh,\Z).$$
The definition we use allows us to give better bounds on
resolution of singularities required.

\subsection{Hodge and de Rham cohomology}
Consider the sheaf of differentials $\Omega^n$ on $\Sm^d/k$ and its pull-back 
$\rho^*\Omega^n$ to $(\Sch^d/k)_\eh$. 
Note that $\rho^*\Omega^n$ does not agree with the sheaf of differentials 
for non-smooth schemes. In this section  we study the Hodge cohomology groups
$H^i(X_\eh,\rho^*\Omega^n)$
and de Rham cohomology groups $H^i_{DR}(X_\eh)=H^i(X_\eh,\rho^*\Omega^\cdot)$
for the \eh-topology. We will need Hodge cohomology groups
with compact support in the formula for $\zeta$-values below.
On the other hand, the \eh-version of
de Rham cohomology generalizes
Hartshorne's de Rham cohomology \cite{hartshorne} for fields
of characteristic $0$.
We start with the following analog of Theorem \ref{sameasbefore}.

\begin{theorem}\label{sameasbefore1}
If $X\in Sm^d/k$ and condition $R(d)$ holds, then
\begin{align*}
H^i(X,\Omega_X^n)&\cong H^i(X_\eh,\rho^*\Omega^n);\\
H^i_{DR}(X)&\cong H^i_{DR}(X_\eh).
\end{align*}
\end{theorem}

\proof
Because $\Omega^n$ is quasi-coherent,
$H^i(X,\Omega^n)\isom H^i(X_\et,\Omega^n)$.
Now the proof works exactly as the proof of
Theorem \ref{sameasbefore}, using the Hodge
cohomology analog of Proposition \ref{blowuples}:
If $X'$ is the blow-up of the smooth scheme $X$
along the smooth center $Z$, then there is a long exact
sequence \cite[IV Thm. 1.2.1]{grosthesis}
$$ \cdots \to H^i(X,\Omega^n)\to  H^i(Z,\Omega^n)\oplus
H^i(X',\Omega^n) \to  H^i(Z',\Omega^n)\to \cdots .$$
The statement for de Rham-cohomology follows from the
spectral sequence from Hodge to de Rham cohomology.
\proofend

\begin{corollary}
Assume $R(d)$ and let $X\in Sch^d/k$.
Then the $k$-vector spaces $H^i_c(X_\eh,\rho^*\Omega^n)$ and
$H^i_{DR}(X_\eh)$ are finite dimensional.
\end{corollary}

\proof
This follows from the smooth, projective case by the same argument as Corollary
\ref{finite}.
\proofend

\begin{proposition} Assume $R(d+r)$ and let $X\in Sch^d/k$.

a) (Projective bundle formula) There are canonical isomorphisms
$$
\bigoplus_{j=0}^r H^{i-j}_c(X_\eh,\rho^*\Omega^{n-j})\cong
H^i_c(({\mathbb P}^r_X)_\eh,\rho^*\Omega^n).
$$

b) (Affine bundle formula) There are canonical isomorphisms
$$H^i_c(({\mathbb A}^r_X)_\eh,\rho^*\Omega^n)
\cong H^{i-r}_c(X_\eh,\rho^*\Omega^{n-r}).$$
\end{proposition}

\proof
a) By the usual method we reduce to the smooth and proper case.
In this case, the isomorphism is given by $\sum p_r^*\cup \xi^j$,
where $p_r:{\mathbb P}^r_X\to X$ is the projection and $\xi$ the image of
${\cal O}(1)$ under the map
$\Pic {\mathbb P}^r_X \cong H^1(({\mathbb P}^r_X)_\et,{\mathbb G}_m)
\stackrel{dlog}{\longrightarrow}H^1(({\mathbb P}^r_X)_\et,\Omega^1)$,
see \cite[\S 6.3]{grothendieckdix}.

b) Follows formally from a) as in Proposition \ref{projb}.
\proofend

Recall the definition of algebraic de Rham cohomology of Hartshorne.
Given a scheme $X\in \Sch/k$, we can use a \v Cech-covering argument and
assume that there exists a closed immersion of $X$ into a smooth scheme
$W\in \Sm/k$. Then the de Rham cohomology of $X$ is defined
as the hypercohomology of the formal completion of the de Rham complex
$\Omega^\cdot_W$ of $W$ along $X$ \cite[II \S 1]{hartshorne}
$$H^i_{DR}(X)=H^q(\hat W,\hat \Omega^\cdot_W).$$
If $k$ is of characteristic $0$, then this is independent of
the choice of $W$ by \cite[Thm. 1.4]{hartshorne}.

\begin{theorem}
If $X$ is a scheme of characteristic zero, then
$$ H^i_{DR}(X_\eh)\cong H^i_{DR}(X).$$
\end{theorem}

\proof
We proceed by induction on the dimension of $X$.
By \eqref{lex} and \cite[Prop. 4.1]{hartshorne}, both sides
admit a Mayer-Vietoris sequence for a closed cover. Hence
by induction on the number of irreducible components we can
reduce to the case that $X$ is integral. If $X$ is smooth, then both
sides agree with the usual de Rham-cohomology of $X$
by Theorem \ref{sameasbefore1}. In general,
by resolution of singularites,
we can find a blow-up square \eqref{blowup} with $X'$ smooth.
Now we can compare the long exact sequence \eqref{lex}
to the corresponding long exact sequence for de Rham cohomology
\cite[Theorem 4.4]{hartshorne} to complete the proof.
\proofend

\section{Arithmetic cohomology}
From now on we fix a finite field $\F_q$, and denote by
$R(d)$ the existence of resolution of singularities for
schemes over the algebraic closure of $\F_q$. A detailed version of the
following discussion can be found in \cite{ichweil}.

Let $G\subseteq \Gal(\bar \F_q/\F_q)$ be the Weil group,
i.e. the free abelian group of
rank $1$ generated by the Frobenius endomorphism $\varphi$.

\begin{definition}
A Weil-eh-sheaf $\f$ on $\Sch/\F_q$ is an eh-sheaf on $\Sch/\bar \F_q$
together with an action over the Frobenius endomorphism,
$\phi:\f\to \varphi_*\f$, i.e. a compatible family of
isomorphisms $\phi_S:\f(S)\to \f(S\times_{\bar\F_q,\varphi}\bar\F_q)$
for every scheme $S/\bar \F_q$.
We write $(\Sch/\F_q)_\Wh^\sim$ for the topos of Weil-eh-sheaves on
$\Sch/\F_q$.
\end{definition}

Every eh-sheaf $\f$ of $\Sch/\F_q$ gives rise to a Weil-eh-sheaf on
$\Sch/\F_q$ by pulling back $\f$ along $\Spec \bar\F_q\to \Spec \F_q$,
and restricting the resulting Galois action to the Weil group.
Conversely, there is a push-forward map from
Weil-eh-sheaves on $\Sch/\F_q$ to $\eh$-sheaves on $\Sch/\F_q$,
giving a morphism of topoi
$\gamma:(\Sch/\F_q)_\Wh^\sim \to (\Sch/\F_q)_\eh^\sim$.

For a Weil-eh-sheaf $\f$, we define Weil-eh-cohomology $H^i(X_\Wh ,\f)$
as the derived functor of $\f\mapsto \f(\bar X)^G$.
Similarly, we define Weil-eh-cohomology with compact
support $H^i_c(X_\Wh ,\f)$ of the sheaf $\f$ as the cohomology of the complex
$R\Gamma_GR\Gamma_c(\bar X_\eh,\f)$.

Let $e\in H^1((\F_q)_\Wh,\Z)\cong \Z$ be a generator. Since
$e^2\in H^2((\F_q)_\Wh,\Z)=0$, the sequence
\begin{equation}\label{cupeexact}
 \cdots \to H^{i-1}_c(X_\Wh ,\f^\cdot)\stackrel{\cup e}{\longrightarrow}
H^{i}_c(X_\Wh ,\f^\cdot) \stackrel{\cup e}{\longrightarrow}
H^{i+1}_c(X_\Wh ,\f^\cdot) \to \cdots
\end{equation}
is a complex. The results of \cite{ichweil} carry over to the present situation:

\begin{theorem}\label{weilI}
Let $X\in \Sch/\F_q$ and let $\f^\cdot$ be a complex of \eh-sheaves.

a) There are long exact sequences
$$\to H^i_c(X_\eh,\f^\cdot) \to H^i_c(X_\Wh ,\f^\cdot)\to
H^{i-1}_c(X_\eh,\f^\cdot)\otimes\Q \to H^{i+1}_c(X_\eh,\f^\cdot) \to $$

b) If the cohomology sheaves of $\f^\cdot$ are torsion, then
$$H^i_c(X_\eh,\f^\cdot) \cong  H^i_c(X_\Wh,\f^\cdot)$$

c) If the cohomology sheaves of $\f^\cdot$ are uniquely divisible, then
$$H^i_c(X_\Wh ,\f^\cdot) \cong
H^i_c(X_\eh,\f^\cdot)\oplus H^{i-1}_c(X_\eh,\f^\cdot),$$
and cup product with $e$ is given by the matrix
${0\,0\choose 1\, 0}$. In particular, the sequence \eqref{cupeexact}
is exact.
\end{theorem}

Since $G$ has cohomological dimension $1$, the Leray spectral sequence for
composition of functors breaks up into short exact sequences
\begin{equation}\label{lerayses}
0\to H^{i-1}_c(\bar X_\eh,\f)_G \to H^i_c(X_\Wh ,\f)
\to H^i_c(\bar X_\eh,\f)^G\to 0 .
\end{equation}
Because the Leray spectral sequence is multiplicative,
this implies that we have a commutative diagram
\begin{equation}\label{square}
\begin{CD}
H^i_c(X_\Wh ,\f)@>surj. >> H^i_c(\bar X_\eh,\f)^G\\
@V\cup eVV @VVV \\
H^{i+1}_c(X_\Wh ,\f) @<inj. << H^i_c(\bar X_\eh,\f)_G,
\end{CD}
\end{equation}
where the right vertical map is induced by the identity map.

\begin{corollary}\label{ss}
The Galois-modules $H^i_c(\bar X_\eh,\Q(n))$ are semi-simple
at the eigenvalue $1$.
\end{corollary}

\proof
Semi-simplicity at the eigenvalue $1$ is equivalent to the canonical
map
$H^i_c(\bar X_\eh,\Q(n))^G\stackrel{\alpha_i}{\longrightarrow}
 H^i_c(\bar X_\eh,\Q(n))_G$
being an isomorphism. Using \eqref{square} and \eqref{lerayses}
we get the diagram
$$\begin{CD}
0@>>> \im \alpha_{i-1}@= \im(-\cup e) @>>> 0\\
@III @VVV @VVV @VVV \\
0@>>> H^{i-1}_c(\bar X_\eh,\Q(n))_G @>>> \ker(-\cup e)
@>>> \ker \alpha_i @>>> 0.
\end{CD}$$
By Theorem \ref{weilI} c), $\ker(-\cup e)=\im(-\cup e)$,
and the corollary follows.
\proofend

\begin{definition}
For $n\in \Z$, we define arithmetic cohomology with compact support
as Weil-eh-cohomology of the motivic complex,
$$R\Gamma_c(X_\ar,\Z(n))=
R\Gamma_GR\Gamma_c(\bar X_\eh,\Z(n)).$$
Arithmetic cohomology with coefficients in $A$ is defined as
the Weil-eh-cohomology of the complex $\Z(n)\otimes A$.
\end{definition}

Recall that the Weil-etale-cohomology groups $H^i_W(X,\Z(n))$
of \cite{ichweil} are defined as the cohomology groups of the complex
$R\Gamma(G,R\Gamma(\bar X_\et,\Z(n)))$.
In order to apply the results of \cite{ichweil}, we record the
following consequence of Theorem \ref{sameasbefore}:

\begin{corollary}\label{altgleich}
If $X$ is smooth and projective, then under $R(\dim X)$ we have
$H^i_W(X,\Z(n))\cong H^i_c(X_\ar,\Z(n))$.
\end{corollary}

\begin{lemma}
If $A$ is torsion free, then there are isomorphisms
\begin{align*}
H^i_c(\bar X_\eh,\Z(n))^G\otimes A&\isom
H^i_c(\bar X_\eh,A(n))^G,\\
H^i_c(\bar X_\eh,\Z(n))_G\otimes A&\isom
H^i_c(\bar X_\eh,A(n))_G.
\end{align*}
\end{lemma}

\proof
By torsion freeness we have $H^i_c(\bar X_\eh,\Z(n))\otimes A\cong
H^i_c(\bar X_\eh,A(n))$. 
The map on coinvariants is an isomorphism, because tensor product
and coinvariants commute. For invariants, we compare the sequence
\eqref{lerayses} with $\Z(n)$ and $A(n)$-coefficients:
$$\begin{CD}
H^{i-1}_c(\bar X_\eh,\Z(n))_G\otimes A @>>>
H^i_c(X_\ar ,\Z(n))\otimes A @> surj.>> H^i_c(\bar X_\eh,\Z(n))^G\otimes A\\
@VVV @| @VVV \\
H^{i-1}_c(\bar X_\eh,A(n))_G@>inj. >>
H^i_c(X_\ar ,A(n)) @>surj.>> H^i_c(\bar X_\eh,A(n))^G.
\end{CD}$$
\proofend

The following corollary is an immediate consequence of Proposition
\ref{projb}.

\begin{corollary} Assume $R(d+r)$, let $X\in Sch^d/\F_q$ and $n\in \Z$.

a) (Projective bundle formula) There are canonical isomorphisms
$$ \bigoplus_{j=0}^r H^{i-2j}_c(X_\ar,\Z(n-j))\cong
H^i_c(({\mathbb P}^r_X)_\ar,\Z(n)).$$

b) (Affine bundle formula) There are canonical isomorphisms
$$ H^{i}_c(({\mathbb A}^r\times X)_\ar ,\Z(n))
\cong H^{i-2r}_c(X_\ar ,\Z(n-r)).$$
\end{corollary}

\begin{proposition}\label{boundedabove}
Let $X\in \Sch^d/\F_q$ and assume $R(d)$.

a) $H^{i}_c(X_\ar ,\Z(n))=0$ for $i>\max\{2d+1,n+d+1\}$.

b) If $r$ is the number of irreducible components of $X$
of maximal dimension $d$, then there is a canonical surjection
$H^{2d+1}_c(X_\ar ,\Z(d))\to \Z^r$.
\end{proposition}

\proof
a) This follows from the smooth projective case \cite{ichweil}
with the argument of Lemma \ref{devissage}.

b) Let $X=Y\cup\bigcup_i X_i$, where $Y$ is the union of irreducible
components of dimension smaller than $d$ and the $X_i$ are the irreducible
components of dimension $d$. Then by \eqref{lex}, a) and induction on the
number of irreducible components, we get a surjection
$H^{2d+1}_c(X_\ar ,\Z(d))\to \oplus H^{2d+1}_c((X_i)_\ar ,\Z(d))$, hence we
can assume that $X$ is irreducible.
If $W$ is a compactification of $X$, then by \eqref{cs} and a) we get a
surjection $H^{2d+1}_c(X_\ar ,\Z(d))\to H^{2d+1}(W_\ar ,\Z(d))$, and
this is compatible with maps between compactifications.
Finally, by $R(d)$ we can assume that there is a blow-up
square \eqref{blowup} with $W'$ smooth and projective. Again by a)
and \eqref{lex} the map $H^{2d+1}(W_\ar ,\Z(d))\to H^{2d+1}(W'_\ar ,\Z(d))$
is surjective, and by \cite{ichweil} and Corollary \ref{altgleich}, 
the latter group is $\Z$.
\proofend

\section{Tate's conjecture, finite generation, and
$l$-adic cohomology}
In order not to confuse $l$-adic cohomology ($l\not=p$) with
cohomology with $\Z(n)\otimes\Z_l$-coefficients, we write
$l$-adic cohomology as
\begin{align*} H^i(X_\et,\hat\Z_l(n))&:=
\lim_r H^i(X_\et,\mu_{l^r}^{\otimes n});\\
H^i(X_\et,\hat\Q_l(n))&:= H^i(X_\et,\hat\Z_l(n))\otimes_{\Z_l}\Q_l ,
\end{align*}
and similarly for cohomology with compact support.
These groups agree with the continuous cohomology groups of Jannsen
\cite{jannsencont}. By Theorem \ref{susvoeequal}, we could define $l$-adic
cohomology with the eh-topology.
Assume $R(\dim X)$ and consider the following natural map
\begin{equation}\label{kahnmap} 
H^i_c(X_\ar,\Z_l(n))\to H^i_c(X_\ar,\Z/l^r(n))
\stackrel{\sim}{\gets}H^i_c(X_\eh,\Z/l^r(n))\cong 
H^i_c(X_\et,\mu_{l^r}^{\otimes n}).
\end{equation}
The isomorphisms come from Lemma \ref{sheafisom} and Theorem \ref{weilI} b).

In the following conjectures, let $X\in \Sch/\F_q$ and $n\in \Z$.

\medskip

\noindent{\bf Conjecture K(X,n)}\;{\it
The groups $H^{i}_c(X_\ar ,\Z(n))$ form an integral model for $l$-adic
cohomology with compact support for all $l\not= p$, i.e. the limit
of the maps \eqref{kahnmap} induces an isomorphism
$$H^i_c(X_\ar ,\Z_l(n)) \stackrel{\sim}{\longrightarrow}
H^i_c(X_\et,\hat \Z_l(n)). $$}

\medskip

\noindent{\bf Conjecture L(X,n)}\;{\it
For all $i$, the groups $H^{i}_c(X_\ar ,\Z(n))$ are finitely generated
abelian groups.}

\smallskip

\begin{proposition}\label{uist} Assume $R(d)$. 

a) Conjecture $L(X,n)$ for $X\in \Sch^d/\F_q$ implies $K(X,n)$,
the finiteness of $H^i_c(X_\ar ,\Z(n))$ for $i\not\in [2n,n+d+1]$,
and the vanishing for $i\not\in [0,2d+1]$.

b) Conjecture $L(X,n)$ for all smooth and projective
$X\in Sch^d/\F_q$ implies $L(X,n)$ for all $X\in Sch^d/\F_q$.
The same statement is true for $K(X,n)$.
\end{proposition}

\proof a) The finite generation of $H^{i}_c(X_\ar ,\Z(n))$ implies by the 
long exact coefficent sequence that 
$$H^{i}_c(X_\ar ,\Z_l(n)) \cong \lim_r  H^{i}_c(X_\ar ,\Z(n))/l^r
\cong \lim_r H^{i}_c(X_\ar ,\Z/l^r(n)).$$
The latter group is isomorphic  to $H^i_c(X_\et,\hat \Z_l(n))$ as in 
\eqref{kahnmap}. The $l$-adic cohomology
groups $H^i_c(X_\et,\hat \Z_l(n))$ are finite or zero outside
the given intervals by \cite[Thm. 3]{kahnetale}.

b) This follows easily with Lemma \ref{devissage}.
\proofend

For a smooth and projective variety $X$, let $CH^n(X)$
and $A_{num}^n(X)$ be the free abelian group on
closed integral subschemes of $X$ of codimension $n$ modulo rational
and numerical equivalence, respectively.

\medskip

\noindent{\bf Conjecture (Tate/Beilinson)}
\label{tatebeil}{\it
For all smooth and projective varieties $X/\F_q$
and all $n\in \Z^{\geq 0}$, rational and numerical equivalence for
algebraic cycles of codimension $n$ on $X$ agree up to
torsion, and the order of the pole of the zeta function $\zeta(X,s)$ at
$s=n$ is equal to the rank of $A_{num}^n(X)$:
$$ \dim CH^n(X)\otimes\Q=\dim A^n_{num}(X)\otimes\Q=\ord_{s=n}\zeta(X,s).$$}
\medskip

By Tate \cite[Thm. 2.9]{tate}, this implies that for all smooth
and projective $X$, the cycle map
$$ CH^n(X)\otimes\Q_l \to
H^{2n}(\bar X_\et,\hat \Q_l(n))^{\Gal(\bar \F_q/\F_q)}$$
is an isomorphism, and that $H^{2n}(\bar X_\et,\hat \Q_l(n))$
is semi-simple at eigenvalue $1$.

\begin{theorem}
If $R(d)$ and the Tate-Beilinson conjecture hold, then $L(X,n)$
holds for all $X\in Sch^d/\F_q$ and all $n$.
\end{theorem}

\proof
In view of Corollary \ref{altgleich}, we get $L(X,n)$ for
all smooth and projective varieties in $\Sch^d/\F_q$
from \cite[Thm. 8.4]{ichweil}. The general case follows by
Proposition \ref{uist} b).
\proofend

The following conjecture can be thought of as the dual of the generalized 
Tate conjecture of Jannsen \cite[Conj. 12.4, 12.6]{jannsenln} for 
arbitrary $X\in \Sch/\F_q$ and $n\in \Z$.
\medskip

\noindent{\bf Conjecture J(X,n)}\;{\it  For all $l\not=p$, the canonical map
\begin{equation}\label{tate1}
H^i_c(X_\eh,\Q_l(n))\to
H^i_c(\bar X_\et,\hat \Q_l(n))^{\Gal(\bar \F_q/\F_q)}
\end{equation}
is an isomorphism, and the Galois-module
$H^i_c(\bar X_\et,\hat \Q_l(n))$ is semi-simple at the eigenvalue $1$.}
\medskip

We have $CH^n(X)\otimes \Q\isom H^{2n}(X_\et,\Q(n))$ because rationally
Zariski and etale motivic cohomology agree. Hence for smooth and 
projective $X$, $J(X,n)$ specializes to Tate's conjecture on the bijectivity 
of the cycle map under resolution of singularities, because then
$H^{2n}(X_\et,\Q(n))\isom H^{2n}(X_\eh,\Q(n))$ by Theorem \ref{sameasbefore}.
By Proposition \ref{counterex} a) below, $J(X,n)$ is wrong if one uses
the etale topology.

\begin{lemma}\label{ladic}
We have $\dim H^i_c(\bar X_\et,\hat \Q_l(n))^G=
\dim H^{i}_c(\bar X_\et,\hat \Q_l(n))_G $, and there is a short exact sequence
of finite dimensional $\Q_l$-vector spaces
$$ 0\to H^{i-1}_c(\bar X_\et,\hat \Q_l(n))_G
\to H^i_c(X_\et,\hat \Q_l(n))\to H^i_c(\bar X_\et,\hat \Q_l(n))^G \to 0.$$
\end{lemma}

\proof
The first statement follows from the exact sequence of finite
dimensional $\Q_l$-vector spaces:
\begin{multline*}
0\to H^i_c(\bar X_\et,\hat \Q_l(n))^G \to
H^i_c(\bar X_\et,\hat \Q_l(n))
\stackrel{1-\varphi}{\longrightarrow}\\
H^i_c(\bar X_\et,\hat \Q_l(n)) \to
H^i_c(\bar X_\et,\hat \Q_l(n))_G \to 0.
\end{multline*}
Taking the inverse limit over $r$ of the short exact sequences of finite
groups
\begin{multline*}
0\to H^i_c(\bar X_\et,\mu_{l^r}^{\otimes n})^G \to
H^i_c(\bar X_\et,\mu_{l^r}^{\otimes n})
\stackrel{1-\varphi}{\longrightarrow}\\
H^i_c(\bar X_\et,\mu_{l^r}^{\otimes n}) \to
H^i_c(\bar X_\et,\mu_{l^r}^{\otimes n})_G \to 0,
\end{multline*}
and comparing with the kernel and cokernel of $1-\varphi$ on
$H^i_c(\bar X_\et,\hat \Q_l(n))$, we get
$\lim \big(H^i_c(\bar X_\et,\mu_{l^r}^{\otimes n})_G\big)\otimes \Q_l
\cong H^i_c(\bar X_\et,\hat \Q_l(n))_G$,
and
$\lim \big(H^i_c(\bar X_\et,\mu_{l^r}^{\otimes n})^G\big)\otimes\Q_l
\cong H^i_c(\bar X_\et,\hat \Q_l(n))^G$.
Taking the inverse limit of
$$ 0\to H^{i-1}_c(\bar X_\et,\mu_{l^r}^{\otimes n})_G
\to H^i_c(X_\et,\mu_{l^r}^{\otimes n})\to
H^i_c(\bar X_\et,\mu_{l^r}^{\otimes n})^G \to 0,$$
the Lemma follows.
\proofend

\begin{theorem}\label{main1} Assume $R(\dim X)$. Then
$K(X,n) \Leftrightarrow J(X,n)$, and $J(X,n)$ for all
smooth and projective $X\in Sch^d/\F_q$ implies $J(X,n)$
for all $X\in Sch^d/\F_q$.
\end{theorem}

The second statement reproves Theorem 12.7 of Jannsen
\cite{jannsenln}.

\proof
Assuming $K(X,n)$, Theorem \ref{weilI} c), and multiplicativity of the
cycle map imply that cup product with the image of $e$ under the canonical map
$\Z\cong H^{1}((\F_q)_\ar ,\Z(0))\to \Q_l\cong H^1((\F_q)_\et,\hat \Q_l(0))$
gives a long exact sequence of $l$-adic cohomology groups
$H^*_c(X_\et,\hat \Q_l(n))$. Using the short exact sequences of Lemma
\ref{ladic}, the argument of Corollary \ref{ss} shows
semi-simplicity of $H^i_c(\bar X_\et,\hat \Q_l(n))$ at the eigenvalue
$1$ for all $i$.

We now prove by induction on $i$ that the map \eqref{tate1} is an
isomorphism. By Corollary \ref{ss} and semi-simplicity, we have
a commutative diagram
\begin{equation}\label{dh1}\begin{CD}
H^{i-1}_c(X_\eh,\Q_l(n))@>>>
H^{i-1}_c(\bar X_\eh,\Q_l(n))^G @>\sim >> H^{i-1}_c(\bar X_\eh,\Q_l(n))_G \\
@III @VVV @VVV \\
@EEE H^{i-1}_c(\bar X_\et,\hat \Q_l(n))^G @>\sim >>
H^{i-1}_c(\bar X_\et,\hat \Q_l(n))_G.
\end{CD}\end{equation}
Since $H^{i-1}_c(\bar X_\eh,\Q_l(n))^G\cong
H^{i-1}_c(\bar X_\eh,\Q_l(n))^{\Gal(\bar \F_q/\F_q)}$, and
higher Galois cohomology is torsion, a Hochschild-Serre spectral
sequence argument shows that the upper left map is an isomorphism. Hence the
induction hypothesis implies that the right hand map
is an isomorphism. By \eqref{lerayses} and Lemma \ref{ladic},
we have a commutative diagram with exact rows
\begin{equation}\label{dh2}\begin{CD}
H^{i-1}_c(\bar X_\eh,\Q_l(n))_G@>inj.>> H^{i}_c(X_\ar,\Q_l(n))
@>surj.>> H^i_c(\bar X_\eh,\Q_l(n))^G \\
@VVV @V\cong VV @VVV \\
H^{i-1}_c(\bar X_\et,\hat \Q_l(n))_G@>inj.>> H^i_c(X_\et,\hat \Q_l(n))
@>surj.>> H^i_c(\bar X_\et,\hat \Q_l(n))^G.
\end{CD}\end{equation}
Because the left hand map is an isomorphism, so is the right hand map, which
concludes the induction step.

To prove the converse implication, note that by Lemma \ref{sheafisom}, it
suffices to show $K(X,n)$ after tensoring with $\Q_l$.
By hypothesis, the left hand vertical map, hence
the right hand vertical map in \eqref{dh1} is an isomorphism
for all $i$, and we conclude that the middle vertical map in \eqref{dh2}
is an isomorphism. The final statement follows from Proposition \ref{uist}. 
\proofend

\section{Values of zeta-functions}
We give a conjectural formula for values of zeta functions in terms of
arithmetic cohomology, inspired by Lichtenbaum \cite{licht}.
Fix $X\in \Sch/\F_q$, let $\zeta(X,s)$ be the zeta function
of $X$, $\chi(H^{*}_c(X_\ar ,\Z(n)),e)$ the Euler-characteristic of the
complex
$$\cdots \to H^{i-1}_c(X_\ar ,\Z(n))\stackrel{\cup e}{\longrightarrow}
H^{i}_c(X_\ar ,\Z(n))\stackrel{\cup e}{\longrightarrow}
H^{i+1}_c(X_\ar ,\Z(n))\to\cdots,$$
and $$\chi(n)= \sum_{0\leq i\leq n\atop j}(-1)^{i+j}(n-i)\cdot 
\dim H^j_c(X_\eh,\rho^*\Omega^i).$$

\medskip
\noindent{\bf Conjecture Z(X,n)}\;{\it
The alternating sum $\sum_i (-1)^i\rk H^{i}_c(X_\ar ,\Z(n))$ is zero,
the weighted alternating sum of the ranks equals the order of the
zeta-function
$$\sum_i (-1)^i i\cdot \rk H^{i}_c(X_\ar ,\Z(n))=
\ord_{s=n}\zeta(X,s)=:\rho_n,$$
and for $s\mapsto n$,
\begin{equation}\label{zetavalue1}
\zeta(X,s)\sim\pm
(1-q^{n-s})^{\rho_n}\cdot \chi(H^{*}_c(X_\ar ,\Z(n)),e)
\cdot q^{\chi(n)}.
\end{equation}}
\medskip

Note that the finiteness of the sum, and the finiteness of the sum
defining the Euler characteristic is part of the conjecture.

\begin{theorem}\label{zetavalue}
Under $R(d)$, $L(X,n)$ for all $X\in Sch^d/\F_q$ implies $Z(X,n)$ for all
$X\in Sch^d/\F_q$.
\end{theorem}

\proof By $L(X,n)$ and Theorem \ref{weilI} c), the alternating sum
of the ranks of the groups $H^{i}_c(X_\ar ,\Z(n))$ equals
$$\sum_i (-1)^i\big(\rk H^i_c(X_\eh,\Z(n))+
\rk H^{i-1}_c(X_\eh,\Z(n))\big)=0.$$
From Lemma \ref{ladic} and Proposition \ref{uist} a), we get
\begin{multline*}
\sum_i (-1)^i i \cdot \rk H^{i}_c(X_\ar ,\Z(n))\\
=\sum_i (-1)^i i\cdot\big( \dim H^{i-1}_c(\bar X_\et,\hat \Q_l(n))_G
+\dim H^i_c(\bar X_\et,\hat \Q_l(n))^G \big)\\
=-\sum_i(-1)^i \dim H^i_c(\bar X_\et,\hat \Q_l(n))^G
= -\sum_i(-1)^i \dim H^i_c(\bar X_\et,\hat \Q_l)^{\varphi=q^n}
\end{multline*}
By Grothendieck's formula for $\zeta(X,s)$, the latter
agrees with $\ord_{s=n}\zeta(X,s)$.

For smooth and projective varieties, the formula \eqref{zetavalue1}
holds by \cite[Thm. 9.1]{ichweil} and Corollary \ref{altgleich}.
For arbitrary $X$, note that by $L(X,n)$ and Theorem \ref{weilI} c),
the cohomology groups of the complex $(H^{*}(X,\Z(n)),e)$ are finite.
By the argument of Lemma \ref{devissage}, it suffices to show the
following: If $U\subseteq X$
is an open subscheme with complement $Z$, and the formula \eqref{zetavalue1}
holds for two of the three schemes $Z$, $X$ and $U$, then it also holds for 
the third. Consider the double complex
\begin{equation}\begin{CD}\label{zetass}
@A\cup eAA @A\cup eAA @A\cup eAA \\
\cdots\to H^{i+1}_c(U_\ar ,\Z(n))@>>>  H^{i+1}_c(X_\ar ,\Z(n))@>>>
H^{i+1}_c(Z_\ar ,\Z(n))\to\cdots
\\
@A\cup eAA @A\cup eAA @A\cup eAA \\
\cdots \to H^i_c(U_\ar,\Z(n))@>>>  H^i_c(X_\ar,\Z(n))@>>>H^i_c(Z_\ar,\Z(n))
\to \cdots\\
@A\cup eAA @A\cup eAA @A\cup eAA
\end{CD}
\end{equation}
Every anti-diagonal has only finitely many non-zero entries by
Proposition \ref{uist} a). Taking horizontal cohomology, we see that the
double complex is exact.
Taking vertical cohomology, we get the $E_1$-terms of a spectral
sequence whose $E_1$-terms are finite, which converges to zero, and which 
has only finitely many differentials, i.e. $E_r=E_\infty$ for $r>> 0$.
An inspection shows that the equality
$$ \chi(H^{*}_c(X_\ar ,\Z(n)),e)=\chi(H^{*}_c(Z_\ar ,\Z(n)),e)\cdot
\chi(H^{*}_c(U_\ar ,\Z(n)),e)$$
is equivalent to the product of the orders of the $E_1$-terms
on a anti-diagonal being equal for two adjacent anti-diagonals, i.e.
$\prod_i |E^{i,-i}_1|=\prod_i |E^{i+1,-i}_1|$.
But it is easy to see that this property is preserved under differentials, i.e.
$\prod_i |E^{i,-i}_r|=\prod_i |E^{i+1,-i}_r|$
if and only if $\prod_i |E^{i,-i}_{r+1}|=\prod_i |E^{i+1,-i}_{r+1}|$.
Now the claim follows because the spectral sequence converges
to zero, and $E_r=E_\infty$ for $r>>0$, hence both sides equal one
for $r>>0$.

For the $p$-part, it is easy to see that for fixed $i$,
$\sum_j(-1)^j\dim H^j_c(X_\eh,\rho^*\Omega^i)$ is compatible with
the localization sequences \eqref{cs}, thus in view of
Theorem \ref{sameasbefore}, we can deduce the result
from the smooth and projective case \cite[Thm. 0.1]{milnevalues}.
\proofend

\rem
We cannot use de Jong's Theorem on alterations to prove Theorem 
\ref{zetavalue}, because it is not clear how the
formula $Z(X,n)$ behaves under finite etale Galois extensions.
Also, it does not suffice to assume $K(X,n)$ instead of $L(X,n)$. Indeed,
one can construct a diagram of the form \eqref{zetass},
with torsion vertical cohomology groups,
where all vertical cohomology groups for two of the three complexes are zero,
but the vertical cohomology of the third complex is not finitely generated
(because if one considers the spectral sequence to the double complex
\eqref{zetass}, one can have a differential $d_3$ that is non-trivial for
infinitely many primes).

\section{Examples}

\begin{theorem}
a) If $R(d)$ holds, then $L(X,n)$ holds for
all $X\in \Sch^d/\F_q$ and $n\leq 0$.

b) If $\dim X\leq 1$, then $L(X,n)$ holds for all $n\in \Z$.

c) Let $X$ be a surface for which every irreducible component is
birationally equivalent to a surface satisfying Tate's conjecture. Then
$L(X,n)$ and $Z(X,n)$ holds for $n\leq 1$.
\end{theorem}

\proof
a) For $X$ smooth and projective this is \cite[Prop. 9.2]{ichweil}
and  Corollary \ref{altgleich}.
The general case follows with Lemma \ref{devissage}.

b) The statement is easy for zero-dimensional schemes, so it suffices
by Lemma \ref{devissage} to show the statement for a smooth and
projective curve. But then the result is the combination of \cite[Prop. 9.4]{ichweil},
\cite[Thm. 8.4]{ichweil} and Corollary \ref{altgleich}.

c) By the curve case and Lemma \ref{devissage}, we can assume that $X$ is
smooth and projective. In this case, we apply \cite[Thm. 9.3]{ichweil} for
$n=1$, and a) for $n\leq 0$.
\proofend

\noindent{\bf Example 1.} (Zero-dimensional schemes) Since arithmetic 
cohomology groups and zeta functions are invariant under nilpotent
extensions and compatible with coproducts, it suffices to consider the 
case $X=\Spec \F_{q^r}$. The zeta-function is
$$\zeta(\F_{q^r},s)=\zeta(\F_q,rs)=\frac{1}{1-q^{-rs}}.$$
Let $w_{rn}=|\Q/\Z(n)^{\Gal(\bar \F_q/\F_{q^r})}|=
q^{|rn|}-1$ if $n\not=0$, and $w_0=1$. Then
\begin{equation}\label{field}
H^{i}_c((\F_{q^r})_\ar ,\Z(n))=
\begin{cases}
\Z &n=0, i=0,1;\\
\Z/w_{rn} &n\not=0, i=1 ;\\
0 & \text{otherwise}.
\end{cases}
\end{equation}
This is clear for $n=0$. For $n\not=0$, we use that
$H^{i}_c((\F_{q^r})_\ar ,\Z(n))=H^{i-1}((\F_{q^r})_\et,\Q/\Z(n))$,
because cohomology with rational coefficients vanishes.
For the $p$-part, one checks easily that $\chi(n)=rn$ for $n\geq 0$ and
$\chi(n)=0$ for $n<0$. 
For $n\not=0$, formula \eqref{zetavalue1} becomes the identity
$\frac{1}{1-q^{-rn}}=\pm (1-q^{n-s})^0 w_n^{-1}\cdot q^{\chi(n)} $. 
For $n=0$, we have
$\chi(H^*_c((\F_{q^r})_\ar ,\Z(0)),e)=\frac{1}{r}$, because the map 
$H^0_c((\F_{q^r})_\ar ,\Z)\to H^1_c((\F_{q^r})_\ar ,\Z)$
is an injection with cokernel $\Z/r$. Indeed,
it is induced by the map $(\Z^r)^G\to (\Z^r)_G$, and the former is
generated by $(1,\ldots,1)$, whereas the latter is generated by
the class of $(1,0,\ldots,0)$. Formula \eqref{zetavalue1} holds
because $\lim_{s\to 0}\frac{1-q^{-rs}}{1-q^{-s}}=r$.

\bigskip

\noindent{\bf Example 2.}
Let $C={\mathbb P^1}/0\sim 1$ be the node, and ${\mathbb P^1}$ be its
normalization. Consider the blow-up square
\begin{equation}
\begin{CD}\label{curve}
\Spec \F_q\amalg \Spec \F_q @>\iota>>{\mathbb P^1} \\
@Vi'VV @VVV \\
\Spec \F_q @>>> C.
\end{CD}
\end{equation}
The corresponding long exact sequence \eqref{lex} breaks up into short exact
sequences
\begin{equation}\label{ujuj}
0\to H^{i-1}_c((\F_p)_\ar ,\Z(n))\to H^i_c(C_\ar,\Z(n))
\to H^i_c((\P^1)_\ar,\Z(n))\to 0,
\end{equation}
because both maps $\iota^*$ and ${i'}^*$ have image the diagonal
of $H^{i}_c((\F_p)_\ar ,\Z(n))\oplus H^i_c((\F_p)_\ar ,\Z(n))$.
Using the projective bundle formula and Example 1, we get
\begin{equation}\label{qqqq}
\rk H^i_c(C_\ar,\Z(n)) =\begin{cases}
2 & (n,i)=(0,1);\\
1 & (n,i)=(0,0),(0,2),(1,2),(1,3) ;\\
0 & \text{otherwise}.
\end{cases}
\end{equation}
For the weighted alternating sum of the ranks we get
$\rho_1=-1$, and $\rho_n=0$ for $n\not=0$. For the precise value of the
zeta-function, one calculates from \eqref{ujuj}, Example 1 and the
projective bundle formula
$$|H^{i}_c(C_\ar ,\Z(n))_{tor}|=
\begin{cases}
w_n &i=1,2 ;\\
w_{n-1} & i=3;\\
0 &\text{otherwise}.
\end{cases}$$
This gives $\chi(H^{*}_c(C_\ar ,\Z(n)),e)=w_{n-1}^{-1}$. The same
argument shows that $\chi(n)=n-1$ for $n>0$ and $\chi(n)=0$
for $n\leq 0$, hence for $s\mapsto n$
$$\zeta(C,s)=\frac{1}{1-q^{1-s}}=
\pm (1-q^{n-s})^{\rho_n}\cdot q^{\chi(n)}\cdot w_{n-1}^{-1}.$$

\bigskip

\noindent{\bf Example 3.}
We give an example which shows that  all conjectures above are wrong
if statet for the etale topology instead of the eh-topology
(we received help from C.Weibel in constructing this example).
Let $X$ be a normal surface over $\F_q$ of Reid with one singular point
$P$, such that the blow-up $X'$ of $X$ at $P$ is smooth
and has a node $C$ as the exceptional divisor see
\cite[6.6]{weibelsurface}.

\begin{proposition}\label{counterex}
a) We have
\begin{align*}
H^i(X_\et,\Q)&=\begin{cases}
\Q &i=0;\\
0 & i>0;\end{cases}\\
H^i(X_\eh,\Q)&=\begin{cases}
\Q &i=0, 2;\\
0 & i>0.\end{cases}
\end{align*}

b) Let $i:P\to X$, $i':C\to X'$ be the closed embeddings. Then
$H^2_\et(X,\cone(\Q\to i_*\Q))\cong 0$ and 
$H^2_\et(X',\cone(\Q\to i'_*\Q))\cong \Q$, i.e. etale cohomology with 
compact support of $X-P$ depends on the compactification.

c) We have $\Q/\Z\subseteq H^3_W(X,\Z)$. In particular, Weil-etale
motivic cohomology groups are not finitely generated in general, even for
proper schemes.

d) We have
$$H^i(X_\et,\hat \Q_l)
=\begin{cases}
\Q_l &i\leq 3;\\
0 & i>3.\end{cases}$$
In particular, Conjecture $K(X,0)$ is wrong if we use
etale cohomology instead of eh-cohomology.

e) We have
$$\ord_{s=0}\zeta(X,s)=-2=\sum_i (-1)^ii \cdot \rk H^{i}(X_\ar,\Z(0))$$
but $\sum_i (-1)^i i\cdot \rk H^i_W(X,\Z)=-1$.
In particular, Conjecture $Z(X,0)$ is wrong if one uses
the etale topology.
\end{proposition}

\proof a) Let $g:\Spec F\to X$ be the generic point of $X$. Since
$X$ is normal, we have $g_*\Q\cong \Q$, and $R^sg_*\Q=0$ for
$s>0$, because the stalk of $R^sg_*\Q$ at a geometric point $x\in
X$ is $H^s(K_\et,\Q)$, where $K$ is the field of fractions of the
strictly local ring of $X$ at $x$. The latter group is zero
because higher Galois cohomology is torsion. Hence
$H^i(X_\et,\Q)\cong H^i(F_\et,\Q)=0$ for $i>0$.

To calculate the \eh-cohomology of $\Q$, we apply the long exact
sequence \eqref{lex}, using $H^i(X'_\eh,\Q)=H^i(P_\eh,\Q)=0$
for $i>0$, and the analog of \eqref{qqqq}.

b) This follows by a) from the long exact sequence \eqref{cs} and
$H^1(({\F_q})_\et,\Q)=0$ and $H^1(C_\et,\Q)\cong \Q$,
\cite{weibelsurface}.

c) The long exact coefficient sequence for Weil-etale
cohomology gives
\begin{equation}\label{tyuyt}
\cdots \to H^2_W(X,\Q)\to H^2_W(X,\Q/\Z)\to H^3_W(X,\Z)\to
H^3_W(X,\Q) \to \cdots .
\end{equation}
By a) and the analog of Theorem \ref{weilI} c) the extreme terms vanish, and
using the analog  of Theorems \ref{weilI} b) and \ref{susvoeequal}, we get
$$H^2(X_\eh,\Q/\Z)\cong H^2(X_\et,\Q/\Z)\cong H^2_W(X,\Q/\Z)\cong
H^3_W(X,\Z).$$ But by a), the former group contains 
$\Q/\Z=H^2(X_\eh,\Z)\otimes \Q/\Z$.
(This phenomenon does not occur for the Weil-eh-topology, because
the extreme terms in \eqref{tyuyt} do not vanish.)

d) This can be calculated using blow-up sequences by Lemma \ref{ehsheaf}.

e) Counting points, we get for the zeta-function the formula
$$ \zeta(X,s)\cdot \zeta(C,s)\cong \zeta(X',s)\cdot \zeta(P,s),$$
and since $\ord_{s=0}\zeta(\F_q,s)=-1$,
$\ord_{s=0}\zeta(X',s)=-1$, and
$\ord_{s=0}\zeta(C,s)=\ord_{s=0}\zeta({\mathbb A}^1_{\F_q},s)=0$,
we have $\ord_{s=0}\zeta(X,s)=-2$.
But in view of $H^i_W(X,\Q)\cong H^i(X_\et,\Q)\oplus H^{i-1}(X_\et,\Q)$
we have $-1$ for the weighted alternating sum of Weil-etale
cohomology groups, and $-2$ for the weighted alternating
sum of arithmetic cohomology groups.
\proofend


\begin{thebibliography}{99}
\bibitem{sga4}{\sc M.Artin, A.Grothendieck, J.L.Verdier},
Th\'eorie des topos et cohomologie \'etale des sch\'emas. Tome 1:
Th\'eorie des topos. S\'eminaire de G\'eom\'etrie Alg\`ebrique du
Bois-Marie 1963--1964 (SGA 4).
Lecture Notes in Mathematics, Vol. 269. Springer-Verlag,
Berlin-New York, 1972.

\bibitem{dejong} {\sc A.J. de Jong}, Smoothness, semi-stability and 
alterations. Inst. Hautes \'Etudes Sci. Publ. Math. No. {\bf 83} (1996), 
51--93.

\bibitem{friedvoe} {\sc E.Friedlander, V.Voevodsky}, Bivariant cycle
cohomology. Cycles, transfers, and motivic homology theories, 138--187,
Ann. of Math. Stud., 143, Princeton Univ. Press, Princeton, NJ, 2000.

\bibitem{gabber} {\sc O.Gabber}, Sur la torsion dans la cohomologie
$l$-adique d'une variete. C. R. Acad. Sci. Paris Ser. I Math. {\bf 297} (1983),
no. 3, 179--182.

\bibitem{ichtate} {\sc T.Geisser}, Tate's conjecture, Algebraic Cycles,
and Rational K-theory in characteristic $p$. K-theory {\bf 13} (1998),
109-122.

\bibitem{ichweil} {\sc T.Geisser}, Weil-etale cohomology
over finite fields. Math. Ann. {\bf 330} (2004), 665--692.

\bibitem{ichdede} {\sc T.Geisser}, Motivic cohomology over Dedekind rings.
Math. Z. {\bf 248} (2004), 773--794.

\bibitem{marcI} {\sc T.Geisser, M.Levine}, The $p$-part of
$K$-theory of fields in characteristic $p$. Inv. Math. {\bf 139}
(2000), 459--494.

\bibitem{grosthesis} {\sc M.Gros}, Classes de Chern et classes de cycles en
cohomologie de Hodge-Witt logarithmique.
Mem. Soc. Math. France (N.S.) {\bf 21} (1985).

\bibitem{grossuwa} {\sc M.Gros, N.Suwa}, La conjecture de Gersten pour
les faisceaux de Hodge-Witt logarithmique. Duke Math. J. {\bf 57} (1988),
no. 2, 615--628.

\bibitem{grothendieckdix} {\sc A.Grothendieck},
Classes de Chern et repr\'esentations lin\'eaires des groupes discrets.
1968 Dix Expos\'es sur la Cohomologie des Sch\'emas pp. 215--305
North-Holland, Amsterdam.

\bibitem{hartshorne} {\sc R.Hartshorne}, On the de Rham cohomology
of algebraic varieties. Inst. Hautes Etudes Sci. Publ. Math. No. 
{\bf 45} (1975), 5--99.

\bibitem{jannsenln} {\sc U.Jannsen}, Mixed motives and algebraic
$K$-theory. With appendices by S. Bloch and C. Schoen.
Lecture Notes in Mathematics, {\bf 1400}. Springer-Verlag, Berlin, 1990.

\bibitem{jannsencont}  {\sc U.Jannsen}, Continuous etale cohomology.
Math. Ann. {\bf 280} (1988), 207--245.

\bibitem{kahnold} {\sc B.Kahn}, A sheaf-theoretic reformulation of the
Tate conjecture. Preprint of the Institute de Math\'ematique de Jussieu
no {\bf 150}, 1998.

\bibitem{kahnanens} {\sc B.Kahn}, \'Equivalence rationnelle et num\'erique
sur certaines vari\'et\'es de type ab\'elien sur un corps fini.
Annales Ec. Nor. Sup. {\bf 36} (2003), 977--1002.

\bibitem{kahnetale} {\sc B.Kahn}, Some finiteness results for \'etale
cohomology. J. Number Theory {\bf 99} (2003), 57-73.

\bibitem{marcbook} {M.Levine}, Mixed motives. Mathematical Surveys and
Monographs, {\bf 57}. American Mathematical Society, Providence, RI, 1998.

\bibitem{licht} {\sc S.Lichtenbaum}, The Weil-etale topology.
To appear in: Compositio Math.

\bibitem{milnebook} {\sc J.S.Milne}, Etale cohomology. Princeton Math.
Series 33.

\bibitem{milnevalues} {\sc J.S.Milne},
Values of zeta functions of varieties over finite fields.
Amer. J. Math. {\bf 108} (1986), no. 2, 297--360.


\bibitem{sv96} {\sc A.Suslin, V.Voevodsky}, Singular homology of abstract
algebraic varieties. Invent. Math. {\bf 123} (1996), no. 1, 61--94.

\bibitem{sv00} {\sc A.Suslin, V.Voevodsky}, Bloch-Kato conjecture and
motivic cohomology with finite coefficients. The arithmetic and geometry of
algebraic cycles (Banff, AB, 1998), 117--189, NATO Sci. Ser. C Math. Phys.
Sci., 548, Kluwer Acad. Publ., Dordrecht, 2000.

\bibitem{tate} {\sc J.Tate}, Conjectures on algebraic cycles in
$l$-adic cohomology. Motives (Seattle, WA, 1991), 71--83,
Proc. Sympos. Pure Math., {\bf 55}, Part 1,
Amer. Math. Soc., Providence, RI, 1994.

\bibitem{weibelsurface} {\sc C.Weibel}, The negative $K$-theory of normal
surfaces. Duke Math. J. {\bf 108} (2001), no. 1, 1--35.
\end{thebibliography}
\end{document}